\documentclass[12pt]{article}
\usepackage{amsmath, amsthm, amssymb}
\usepackage[dvips]{graphicx}
\usepackage{epsfig}
\usepackage{vector}
\def\eps{\varepsilon}
\font\tencmmib=cmmib10 \skewchar\tencmmib '60
\newfam\cmmibfam
\textfont\cmmibfam=\tencmmib

\def\bbox{\quad\hbox{\vrule \vbox{\hrule \vskip2pt \hbox{\hskip2pt
\vbox{\hsize=1pt}\hskip2pt} \vskip2pt\hrule}\vrule}}
\def\lessim{\ \lower4pt\hbox{$
\buildrel{\displaystyle <}\over\sim$}\ }
\def\gessim{\ \lower4pt\hbox{$\buildrel{\displaystyle >}
\over\sim$}\ }

\def\C{{\cal C}}

\def\A{{\cal A}}

\def\eps{{\varepsilon}}

\def\la{{\Bigl\langle}}
\def\ra{{\Bigr\rangle}}

\def\qed{\hfill\break\rightline{$\bbox$}}
\newcommand{\e}{\mathbb{E}}
\newcommand{\p}{\mathbb{P}}
\newcommand{\Reals}{\mathbb{R}}
\newcommand{\Natural}{\mathbb{N}}

\newtheorem{lemma}{Lemma}
\newtheorem{theorem}{Theorem}
\newtheorem{corollary}{Corollary}

\makeatletter
\@addtoreset{equation}{section}

\makeatother

\font\tencmmib=cmmib10 \skewchar\tencmmib '60
\newfam\cmmibfam
\textfont\cmmibfam=\tencmmib

\def\bbox{\quad\hbox{\vrule \vbox{\hrule \vskip2pt \hbox{\hskip2pt
\vbox{\hsize=1pt}\hskip2pt} \vskip2pt\hrule}\vrule}}
\def\lessim{\ \lower4pt\hbox{$
\buildrel{\displaystyle <}\over\sim$}\ }
\def\gessim{\ \lower4pt\hbox{$\buildrel{\displaystyle >}
\over\sim$}\ }

\def\eps{\varepsilon}

\def\go0{\to 0}

\def\la{\langle}

\def\leftitem#1{\item{\hbox to\parindent{\enspace#1\hfill}}}

\def\qed{\hfill\break\rightline{$\bbox$}}

\def\ra{\rangle}

\def\sg{\sigma}

\def\sg2{\sigma^2}

\def\__{_{\infty}}

\begin{document}

\title{A new representation of the Ghirlanda-Guerra identities with applications.}
\author{Dmitry Panchenko\thanks{Department of Mathematics, Texas A\&M University, email: panchenk@math.tamu.edu. 
Partially supported by NSF grant.}\\
}
\maketitle
\begin{abstract} 
In this paper we obtain a new family of identities for random measures on the unit ball of a separable Hilbert space 
which arise as the asymptotic analogues of the Gibbs measures in the Sherrington-Kirkpatrick and $p$-spin models and 
which are known to satisfy the Ghirlanda-Guerra identities. We give several applications of the new identities to  structural 
results for such measures.

\end{abstract}
\vspace{0.5cm}

\section{Introduction and main results.}

The Gibbs measures in the Sherrington-Kirkpatrick type spin glass models are known to satisfy two asymptotic stability 
properties - the Aizenman-Contucci stochastic stability \cite{AC} and the Ghirlanda-Guerra identities \cite{GG}.
These properties are typically proved by adding a small perturbation term to the Hamiltonian of the model 
(see \cite{Tal-New}, \cite{SG2}), but in some cases can be proved in a strong sense without perturbation
(see \cite{ACh}, \cite{PGGmixed}) using the validity of the Parisi formula and its properties (\cite{Guerra}, \cite{TPF}, \cite{TPM}). 
In this paper, we will prove a new invariance property for asymptotic
Gibbs' measures that satisfy the Ghirlanda-Guerra identities. Even though the idea behind the main result,
Theorem \ref{Th1} below, was originally motivated by the stability property proved in \cite{ACGG} which unified
the Aizenman-Contucci stochastic stability and the Ghirlanda-Guerra identities,  the proof given here is based only 
on the Ghirlanda-Guerra identities. As one application of the new invariance, we will show how one can reduce 
the problem of ultrametricity to a formally weaker condition on the size of non-ultrametric set. As another application, 
we will prove that one can always find infinite sequence of points in the support of the measure such that all their 
overlaps take values in a given set of overlaps of positive probability. This, for example, implies the positivity principle 
of Talagrand (\cite{SG}, \cite{SG2} or \cite{Posit}) and the fact that the support of such measure is not totally bounded. 
For discrete measure that satisfies the Ghirlanda-Guerra identities we will show that the Gram matrix of the points in its
support is weakly exchangeable and independent of the weights. Finally, we will give a couple of examples of new 
identities for the Poisson-Dirichlet distribution $PD(\zeta).$

Let us consider a random probability measure $G$ on the unit ball of a separable Hilbert space $H$.
We will denote by $(\sigma^l)_{l\geq 1}$ an i.i.d. sequence from this measure, by $\la\cdot\ra$
the average with respect to $G^{\otimes \infty}$ and by $\e$ the expectation with respect to the
randomness of $G$. Let $R_{l,l'} = \sigma^l\cdot \sigma^{l'}$ be the scalar product, or overlap,
of $\sigma^l$ and $\sigma^{l'}$. Random measure $G$ is said to satisfy the Ghirlanda-Guerra
identities if for any $n\geq 2,$ any bounded measurable function $f$ of the overlaps $(R_{l,l'})_{l,l'\leq n}$
and any bounded measurable function $\psi$ of one overlap we have
\begin{equation}
\e \la f\psi(R_{1,n+1})\ra = \frac{1}{n}\hspace{0.3mm} \e\la f \ra \hspace{0.3mm} \e\la \psi(R_{1,2})\ra + \frac{1}{n}\sum_{l=2}^{n}\e\la f\psi(R_{1,l})\ra.
\label{GG}
\end{equation}
Another way to express the Ghirlanda-Guerra identities is to say that, conditionally on $(R_{l,l'})_{1\leq l<l' \leq n}$, 
the law of $R_{1,n+1}$ is given by  the mixture 
\begin{equation}
\frac{1}{n} \hspace{0.3mm}\mu + \frac{1}{n}\hspace{0.3mm} \sum_{l=2}^n \delta_{R_{1,l}}
\label{GGgen}
\end{equation}
where $\mu$ is the law of $R_{1,2}$ under the measure $\e G^{\otimes 2}$.
Given $n\geq 1$, consider $n$ bounded measurable functions 
$f_1,\ldots, f_n: \Reals^{n(n-1)/2}\times \Reals\to\Reals$
and define
\begin{equation}
F(\sigma,\sigma^1,\ldots,\sigma^n) = f_1(R^n, \sigma\cdot\sigma^1)+\ldots+f_n(R^n, \sigma\cdot\sigma^n)
\label{F1}
\end{equation}
where we denoted $R^n = (R_{l,l'})_{l<l'\leq n}.$ For $1\leq l\leq n$ we define
\begin{equation}
F_l(\sigma,\sigma^1,\ldots,\sigma^n) = F(\sigma,\sigma^1,\ldots,\sigma^n)
 - f_l(R^n, \sigma\cdot\sigma^l)+ \int \!\! f_l(R^n,x) \,d\mu(x)
\label{F2}
\end{equation}
and for $l\geq n+1$ we define
\begin{equation}
F_l(\sigma,\sigma^1,\ldots,\sigma^n) =F(\sigma,\sigma^1,\ldots,\sigma^n).
\label{F3}
\end{equation}
This definition (\ref{F3}) for $l\geq n+1$ will not be used in the statement but will appear in the proof of our main result.
\begin{theorem}\label{Th1}
Suppose (\ref{GG}) holds and let $\Phi$ be a bounded measurable function of $R^n.$ Then
\begin{equation}
\e\la\Phi\ra =
\e\Bigl\la
\frac{\Phi \exp \sum_{l=1}^{n} F_l(\sigma^l,\sigma^1,\ldots,\sigma^n)}
{\la\exp F(\sigma,\sigma^1,\ldots,\sigma^n)\ra_{\mathunderscore}^n}
\Bigr\ra,
\label{main}
\end{equation}
where the average $\la\cdot\ra_{\mathunderscore}$ in the denominator is in $\sigma$ only for fixed  $\sigma^1,\ldots, \sigma^n$
and the outside average of the ratio is in $\sigma^1,\ldots, \sigma^n$.
\end{theorem}\noindent
When $n=1,$ it is understood that $\Phi$ is a constant. Notice that one can easily recover the original Ghirlanda-Guerra 
identities from (\ref{main}) by taking $f_1 = t\psi$ and $f_2=\ldots = f_n=0$ and computing the derivative at $t=0.$
One can generalize Theorem \ref{Th1} by iteratively applying the principle behind its proof to the new coordinates
and we will describe this generalization in Theorem \ref{Th1prime} in Section \ref{Sec2}. 
Classical form of the Ghirlanda-Guerra identities (\ref{GG}) has been used successfully to prove several results
about the structure of random measures that satisfy these identities (see e.g. \cite{ParT} and Chapter $12$ in \cite{SG2}) and  
we will use the information contained in the new representation (\ref{main}) to prove some new structural results about 
these measures. The ultimate goal would be to show that such measures must be 
ultrametric,  generalizing the results in \cite{PGG}, \cite{PGG2} and \cite{Tal-New} (inspired by \cite{AA}) and
our next result makes a small step in this direction. Measure $G$ is said  to be ultrametric  if for any $q$,
\begin{equation}
G^{\otimes 3}\bigl(\bigl\{(\sigma^1,\sigma^2,\sigma^3) : R_{1,2}\geq q, R_{1,3}\geq q, R_{2,3}< q\bigr\}\bigr)=0.
\label{ultra}
\end{equation} 
Another way to express (\ref{ultra}) is to say that for any two points $\sigma^1,\sigma^2$ sampled from $G$ 
such that  $R_{1,2}=\sigma^1\cdot \sigma^2 \geq q$ the sets
\begin{align}
A(\sigma^1,\sigma^2) &= \bigl\{\sigma : \sigma\cdot \sigma^1 \geq q,\, \sigma\cdot \sigma^2 < q\bigr\},
\nonumber
\\
A'(\sigma^1,\sigma^2) &= \bigl\{\sigma : \sigma\cdot \sigma^1 < q,\, \sigma\cdot \sigma^2 \geq q\bigr\}
\label{AAprime}
\end{align}
have measure zero. The first application of the identities of Theorem \ref{Th1} is the following result 
which says that if one can show that the measure of $A(\sigma^1,\sigma^2)$ is not too big in some sense 
then (\ref{ultra}) holds with probability one.
\begin{theorem}\label{Prop1}
Suppose that (\ref{GG}) holds. If
\begin{equation}
\e\Bigl\la \frac{I(R_{1,2} \geq q)}{(1-G(A(\sigma^1,\sigma^2)))^2} \Bigr\ra <\infty
\label{integ}
\end{equation}
then (\ref{ultra}) holds with probability one.
\end{theorem}\noindent
If the measure $G$ is ultrametric then $G(A(\sigma^1,\sigma^2))=0$ and the above integral is equal to
$\mu([q,1]) = \e\la I(R_{1,2}\geq q)\ra$. Of course, in a given model it might be just as difficult to prove (\ref{integ}) 
as to prove ultrametricity but formally this condition is weaker and Theorem \ref{Prop1} illustrates 
a new possible direction opened by Theorem \ref{Th1}. Let us now define
\begin{equation}
\bar{F}(\sigma,\sigma^1,\ldots,\sigma^n)
=
F(\sigma,\sigma^1,\ldots,\sigma^n)
-
\frac{1}{n} \sum_{l=1}^{n} F_l(\sigma^l,\sigma^1,\ldots,\sigma^n).
\label{H}
\end{equation}
Another consequence of Theorem \ref{Th1} is  the following.
\begin{theorem}\label{Prop2}
Under (\ref{GG}), for any $n\geq 1,$ with probability one over the choice of measure $G$
and for  $G^{\otimes n}$-almost all  $(\sigma^1,\ldots,\sigma^{n})$, the measures of the sets 
$\{\sigma : \bar{F}(\sigma,\sigma^1,\ldots,\sigma^n) > 0\}$ and 
$\{\sigma : \bar{F}(\sigma,\sigma^1,\ldots,\sigma^n) < 0\}$
are either both positive or both equal to zero at the same time. In particular, the inequalities 
$\bar{F} \geq 0$ and $\bar{F} \leq 0$ hold on the sets of positive measure $G.$ 
\end{theorem}\noindent
Let us give several examples of application of Theorem \ref{Prop2}. First, let us apply it to the case of
$n=2$ and the functions $f_1(R^2,x) = I(x \in C)$ and $f_2(R^2,x) = -I(x \in C)$ for a given set $C$. 
If we denote
\begin{align}
C(\sigma^1,\sigma^2) &= \bigl\{\sigma : \sigma\cdot \sigma^1 \in C,\, \sigma\cdot \sigma^2 \not\in C\bigr\},
\nonumber
\\
C'(\sigma^1,\sigma^2) &= \bigl\{\sigma : \sigma\cdot \sigma^1 \not\in C,\, \sigma\cdot \sigma^2 \in C\bigr\}
\end{align}
then (\ref{H}) can be written in this case as
$
\bar{F}(\sigma,\sigma^1,\sigma^2) = I(\sigma\in C(\sigma^1,\sigma^2)) - I(\sigma \in C'(\sigma^1,\sigma^2))
$
and Theorem \ref{Prop2} implies that the measures of the sets $C(\sigma^1,\sigma^2)$ and $C'(\sigma^1,\sigma^2)$ 
are either both positive or both equal to zero at the same time.  In particular, the measures $G(A(\sigma^1,\sigma^2)),$ 
$G(A'(\sigma^1,\sigma^2))$  of the sets in (\ref{AAprime}) are either both positive or both equal to zero. 

Another application of Theorem \ref{Prop2} is as follows. 
Let us take sets $B_l$ for $l\leq n$ such that  at least one $\mu(B_l)>0$  and let us make the choice of functions $f_l(R^n,x)=I(x\not\in B_l).$
Consider any sample $\sigma^1,\ldots,\sigma^n$ from the measure $G$ such that 
$R_{l,l'} \in B_l\cap B_{l'}$ for all $l\not = l'$. Then $f_l(R^n,R_{l,l'}) = 0$ for all $l\not = l'$ and the condition
$\bar{F}\leq 0$ becomes
$$
\sum_{l\leq n} I(\sigma\cdot\sigma^l \not\in B_l) \leq \frac{1}{n}\sum_{l\leq n} \mu((B_l)^c).
$$
If at least one $\mu(B_l)>0$, the right hand side is strictly less than one and the condition is satisfied only
if all $\sigma\cdot\sigma^l \in B_l.$ Thus, we get the following.
\begin{corollary}\label{Cor1}
Suppose at least one $\mu(B_l)>0.$ Under (\ref{GG}), for $G^{\otimes n}$-almost all  $(\sigma^1,\ldots,\sigma^n)$,
if $R_{l,l'} \in B_l\cap B_{l'}$ for all $l\not = l'$ then on a set of positive measure $G$,
$\sigma\cdot\sigma^l \in B_l$ for all $l\leq n.$ 
\end{corollary}\noindent
If $B$ is a set with $\mu(B)>0$ then using Corollary \ref{Cor1} inductively on $n$ we can find
an infinite sequence $(\sigma^l)$ in the support of $G$ such that all overlaps $R_{l,l'}\in B$.
It is known (Theorem 2 in \cite{PGG}) that if the measure $G$ satisfies the Ghirlanda-Guerra identities
and if $q^*$ is the supremum of the support of the distribution of the overlap $R_{1,2}$ under $\e G^{\otimes 2}$
then with probability one $G$ is concentrated on the sphere of radius $\sqrt{q^*}$.
On the other hand, if $\e \la I(R_{1,2}\leq q)\ra >0$ for some $q<q^*$ then taking $B= [-1, q]$
we can find infinite sequence $(\sigma^l)$ in the support of $G$ such that all overlaps $R_{l,l'}\leq q$. 
This means that if the distribution of $R_{1,2}$ is not concentrated on one point $q^*$ then the support of $G$
on the sphere of radius $\sqrt{q^*}$ contains infinitely many points at some positive distance from each other,
which implies the following.
\begin{corollary} 
Under (\ref{GG}), if $\mu(\{q^*\})<1$ then with probability one the support of $G$ is not totally bounded. 
\end{corollary}\noindent
In a related result, it was proved in  \cite{ACh} that under the Aizenman-Contucci stochastic stability 
the support of the measure is not finite dimensional. As another application, using Corollary \ref{Cor1} 
inductively as above with the choice of $B=[-1,-\eps]$ we can recover the positivity principle of Talagrand 
(see \cite{SG},\cite{SG2}  or \cite{Posit}) which states that under (\ref{GG}), we must have $\mu([-1,0))=0.$ 
Indeed, if $\mu([-1, - \eps])>0$ then we can find an infinite sequence $(\sigma^l)$ in the support of $G$ 
such that all overlaps $R_{l,l'}\leq -\eps$ which is impossible since
\begin{equation}
0\leq \|\sum_{l\leq n} \sigma^l \|^2 = \sum_{l,l'\leq n} R_{l,l'}\leq n- n(n-1)\eps<0
\label{posiN}
\end{equation}
for $n$ large enough. This argument gives us an opportunity to emphasize the strength of Corollary \ref{Cor1} and 
Theorem \ref{Th1} by comparing it with the following very elementary proof of positivity principle pointed out to
the author by Michel Talagrand. Given a set $B,$ let us define the event  
$$B_n =\{R_{l,l'} \in B \mbox{ for all } l\not = l' \leq n\}$$ and notice that 
\begin{equation}
I_{B_{n+1}} \geq I_{B_n} - \sum_{l\leq n} I_{B_n} I(R_{l,n+1}\not\in B).
\label{almostUltra}
\end{equation}
The Ghirlanda-Guerra identities (\ref{GG}) imply that 
$$
\e\la I_{B_n} I(R_{l,n+1}\not\in B) \ra = \frac{\mu(B^c)}{n} \e\la I_{B_n} \ra
$$ 
and using (\ref{almostUltra}) we get $\e\la I_{B_{n+1}} \ra \geq \mu(B) \e\la I_{B_n} \ra$ and, by induction, 
$\e\la I_{B_{n+1}} \ra \geq \mu(B)^{n}$. If $\mu(B)>0$ then with positive probability over the randomness
of $G$, $n$ replicas $\sigma^1,\ldots, \sigma^n$ sampled from $G$ belong to the set $B_n$ with positive 
probability.  Taking $B=[-1,-\eps]$ and appealing to (\ref{posiN}) shows that $\mu([-1,-\eps])$ must be zero,
which proves that $\mu([-1,0))$. On the other hand, in the same situation when $\mu(B)>0,$ 
Corollary \ref{Cor1} gave us the same statement with probability one over the randomness of the measure 
$G$ and, moreover, guaranteed that we can find a vector $(\sigma^1,\ldots,\sigma^n)\in B_n$ with all coordinates 
in the support of $G$ starting from any $\sigma^1$, which is a significantly stronger statement.

Our next application of Theorem \ref{Th1} will be for discrete random measures on the Hilbert space $H$ that satisfy 
the Ghirlanda-Guerra identities.
Suppose that $G = \sum_{l\geq 1} v_l \delta_{\xi_l}$ for some random weights $(v_l)$ and random elements $(\xi_l)$ 
on the unit ball of $H$ and assume that $G$ satisfies (\ref{GG}). We assume that the weights $(v_l)$ are arranged in
non-increasing order and denote by $Q=(\xi_l\cdot \xi_{l'})_{l,l'\geq 1}$ the Gram matrix of scalar products of the points 
in the support of $G$. The matrix $Q$ is called weakly exchangeable if 
\begin{equation}
(\xi_{\rho(l)}\cdot\xi_{\rho(l')}) \stackrel{d}{=} (\xi_l\cdot \xi_{l'})_{l,l'\geq 1}
\label{Exch}
\end{equation}
for any permutation $\rho$ of finitely may indices. We will prove the following. 
\begin{theorem} \label{ThEI}
Under (\ref{GG}), the Gram matrix $Q=(\xi_l\cdot \xi_{l'})_{l,l'\geq 1}$ is weakly exchangeable and independent of the weights $(v_l)$.
\end{theorem} 
\noindent Exchangeability of $Q$ under (\ref{GG}) was proved before in Theorem 3 in \cite{PGG} by way of invariance properties of 
$G$ under some random changes of density and here we will give a different proof as an application of Theorem \ref{Th1}
(in the form of Theorem \ref{Th2} below).  Independence of $(v_l)$ and $Q$ under the assumption (\ref{GG}) is a new result
motivated by Theorem 4.2 in \cite{AA} where it was proved under a different assumptions of robust quasi-stationarity and ergodicity 
of $G$ under a family of random changes of density. The corresponding  argument in \cite{AA} is very subtle and, even though 
the proof of Theorem \ref{ThEI} will require some work, we hope that in the end it is much more intuitive. Of course, Theorem 
\ref{ThEI} was also motivated by the fact that this property holds in the ultrametric case of
the Ruelle probability cascades \cite{Ruelle} (see \cite{Bolthausen}).

Finally, let us mention that (\ref{main}) seems new even for the simplest measures satisfying (\ref{GG}) 
in which case it can be written as a family of identities for the Poisson-Dirichlet distribution
$PD(\zeta)$ and we will state a couple of such identities in Section \ref{Sec3}.

\smallskip
 \textbf{Acknowledgement.} The author would like to thank Michel Talagrand for a number of helpful comments
and, in particular, for pointing out an elementary proof of the positivity principle and the fact that the functions 
$f_l$ in Theorem \ref{Th1} can depend on $R^n$.

\section{Invariance principles.}\label{Sec2}

Before we prove Theorem \ref{Th1}, let us formulate a generalization which is obtained by iterating the same idea in a certain sense.
Consider integers $r\geq 1$ and $1\leq n_1<n_2<\ldots<n_r$ and for $p\leq r$ consider sets
$$I_p = \{n_{p-1}+1,\ldots,n_p\}$$ where $n_0 = 0$. The partition $(I_p)$ represents $r$ groups of 
coordinates from $\{1,\ldots,n\}.$ For $p\leq r,$ let $R^{I_p} = (R_{l,l'})_{l<l', l,l'\in I_p}$ denote the array of overlaps 
of coordinates in $I_p$ and for $l\in I_p$ let us consider functions $f_l(R^{I_p},x)$ of $R^{I_p}$ and $x\in \Reals.$
For simplicity of notations we will sometimes keep the dependence of $f_l$ on $R^{I_p}$ for $l\in I_p$ implicit and simply write $f_l(x)$.
For $p\leq r$, we define
\begin{equation}
F^p(\sigma,\sigma^1,\ldots,\sigma^{n_p}) = f_1(\sigma\cdot\sigma^1)+\ldots+f_{n_p}(\sigma\cdot\sigma^{n_p})
\label{Fp1}
\end{equation}
and then define for $l\in I_p$,
\begin{equation}
F_l^p(\sigma,\sigma^1,\ldots,\sigma^{n_p}) = F^p(\sigma,\sigma^1,\ldots,\sigma^{n_p})
 - f_l(R^{I_p},\sigma\cdot\sigma^l)+ \int \!\! f_l(R^{I_p},x) \,d\mu(x),
\label{Fp2}
\end{equation}
and for $l\not\in I_p$,
\begin{equation}
F_l^p(\sigma,\sigma^1,\ldots,\sigma^{n_p}) =F^p(\sigma,\sigma^1,\ldots,\sigma^{n_p}).
\label{Fp3}
\end{equation}
For $p\leq r,$ let us define
\begin{equation}
Z^p = Z^p(\sigma^1,\ldots,\sigma^{n_p}) =
\frac{ \exp \sum_{l\in I_p} F_l^p(\sigma^l,\sigma^1,\ldots,\sigma^{n_p})}
{\la\exp  F^p(\sigma,\sigma^1,\ldots,\sigma^{n_p})\ra_{\mathunderscore}^{|I_p|}}
\label{Zp}
\end{equation}
where, as in Theorem \ref{Th1}, $\la \cdot\ra_{\mathunderscore}$ is the average in $\sigma$ only.
The following holds.
\begin{theorem}\label{Th1prime}
Suppose (\ref{GG}) holds and $\Phi$ is a bounded measurable function of $R^{I_r}.$ Then
\begin{equation}
\e\la\Phi\ra = \e\la Z^1\ldots Z^r \Phi\ra.
\label{maingen}
\end{equation}
\end{theorem}\noindent
Theorem \ref{Th1} corresponds to the case of the one element partition $(I_p)$, i.e. $r=1$. It is essential that 
the function $\Phi$ in (\ref{maingen}) depends only on the coordinates in the last group $I_r$. Note also that,
using replicas, one can rewrite 
\begin{equation}
Z^1\ldots Z^r = 
\frac{\exp \sum_{p\leq r} \sum_{l\in I_p} F_l^p(\sigma^l,\sigma^1,\ldots,\sigma^{n_p})}
{\la \exp \sum_{p\leq r} \sum_{l\in I_p} F^p(\rho^l,\sigma^1,\ldots,\sigma^{n_p}) \ra_{\mathunderscore}},
\label{Zprod}
\end{equation}
where the average $\la\cdot\ra_{\mathunderscore}$ in the denominator is in the new coordinates $(\rho_l)$.

\smallskip\noindent\textbf{Proof of Theorem \ref{Th1}.}
Without loss of generality, let us assume that $\Phi$ takes values in $[0,1]$ and suppose that
$|f_l|\leq L$ for $1\leq l\leq n$ for some large enough $L.$ For $t\geq 0$ let 
\begin{equation}
\varphi(t) = 
\e\Bigl\la
\frac{\Phi \exp \sum_{l=1}^{n} t F_l(\sigma^l,\sigma^1,\ldots,\sigma^n)}
{\la\exp t F(\sigma,\sigma^1,\ldots,\sigma^n)\ra_{\mathunderscore}^n}
\Bigr\ra.
\end{equation}
We will show that the Ghirlanda-Guerra identities (\ref{GG}) imply that this function is constant,
thus, proving the statement of the theorem, $\varphi(0)=\varphi(1).$ If for $k\geq 1$ we denote
$$
D_{n+k} = \sum_{l=1}^{n+k-1}F_l(\sigma^l,\sigma^1,\ldots,\sigma^n)
-(n+k-1) F_{n+k}(\sigma^{n+k},\sigma^1,\ldots,\sigma^n)
$$
then one can easily compute by induction that (recall (\ref{F3}) and that we average in $\sigma$ 
only in the denominator of (\ref{main}))
$$
\varphi^{(k)}(t) = 
\e\Bigl\la
\frac{\Phi D_{n+1}\ldots D_{n+k} 
\exp \sum_{l=1}^{n+k} t F_l(\sigma^l,\sigma^1,\ldots,\sigma^n)}
{\la\exp t F(\sigma,\sigma^1,\ldots,\sigma^n)\ra_{\mathunderscore}^{n+k}}
\Bigr\ra.
$$
First, let us notice that $\varphi^{(k)}(0)=0.$ Indeed, if we denote $\Phi' = \Phi D_1\ldots D_{n+k-1}$
then $\Phi'$ is the function of the overlaps $(R_{l,l'})_{l,l'\leq n+k-1}$ and 
\begin{eqnarray}
&&
\varphi^{(k)}(0) 
=
\e\Bigl\la \Phi' \Bigl(
\sum_{l=1}^{n+k-1}F_l(\sigma^l,\sigma^1,\ldots,\sigma^n)
-(n+k-1) F_{n+k}(\sigma^{n+k},\sigma^1,\ldots,\sigma^n)
\Bigr)
\Bigr\ra
\label{phizero}
\\
&&
=
\sum_{j=1}^{n}
\e\Bigl\la \Phi' \Bigl(
\sum_{l\not = j}^{n+k-1}f_j(R^n, R_{j,l})
+ \int \! \! f_j (R^n,x)\,d\mu(x)
-(n+k-1) f_{j}(R^n,R_{j,n+k})
\Bigr)
\Bigr\ra
=0
\nonumber
\end{eqnarray}
by the Ghirlanda-Guerra identities in the form of (\ref{GGgen}) applied to each term $j$. Now, since $|F_l| \leq L n$
and $|D_{n+k}| \leq 2L (n+k-1)n$ we get
\begin{eqnarray*}
|\varphi^{(k)}(t)|
&\leq&
\Bigl(\prod_{l=1}^{k} 2L(n+l-1)n\Bigr) \,
\e\Bigl\la
\frac{\Phi 
\exp \sum_{l=1}^{n+k} t F_l(\sigma^l,\sigma^1,\ldots,\sigma^n)}
{\la\exp t F(\sigma,\sigma^1,\ldots,\sigma^n)\ra_{\mathunderscore}^{n+k}}
\Bigr\ra
\\
&=&
\Bigl(\prod_{l=1}^{k} 2L(n+l-1)n\Bigr) \,
\e\Bigl\la
\frac{\Phi 
\exp \sum_{l=1}^{n} t F_l(\sigma^l,\sigma^1,\ldots,\sigma^n)}
{\la\exp t F(\sigma,\sigma^1,\ldots,\sigma^n)\ra_{\mathunderscore}^{n}}
\Bigr\ra
\\
&=&
\prod_{l=1}^{k} (n+l-1)\,(2Ln)^k \, \varphi(t).
\end{eqnarray*}
Consider arbitrary $T>0.$ Again, using that $|F_l| \leq L n$ it is obvious that
$\varphi(t) \leq  e^{2L T n^2}$ for $0\leq t\leq T$ and, therefore,
$$
|\varphi^{(k)}(t)| \leq e^{2L T n^2} \frac{(n+k-1)!}{(n-1)!}\,  (2Ln)^k.
$$
By (\ref{phizero}) and Taylor's expansion
$$
|\varphi(t)-\varphi(0)| \leq \max_{0\leq s\leq t} \frac{|\varphi^{(k)}(s)|}{k!}s^k
\leq
e^{2L T n^2}  \frac{(n+k-1)! }{k! \,(n-1)!} (2Ln t)^k.
$$
Letting $k\to \infty$ we get that $\varphi(t)=\varphi(0)$ for $t<(2Ln)^{-1}.$
Therefore, for any $t_0<(2Ln)^{-1}$ we again have 
$\varphi^{(k)}(t_0)=0$ for all $k\geq 1$ and by Taylor's expansion for $t_0\leq t \leq T,$
$$
|\varphi(t)-\varphi(t_0)| \leq \max_{t_0\leq s\leq t} \frac{|\varphi^{(k)}(s)|}{k!}(t-t_0)^k
\leq
e^{2L T n^2}  \frac{(n+k-1)! }{k! \,(n-1)!} (2Ln (t-t_0))^k.
$$
Letting $k\to\infty$ proves that $\varphi(t) = \varphi(0)$ for $0\leq t< 2(2Ln)^{-1}.$
We can continue in the same fashion to prove this equality for all $0\leq t< T$ and 
note that $T$ was arbitrary.
\qed

\noindent
\textbf{Proof of Theorem \ref{Th1prime}.}
The proof is by induction on $r$. Suppose that (\ref{maingen}) is proved for some $r\geq 1.$
In order to make the induction step, we will separate the coordinates in the last group $I_r$ into
two sets coordinates that will play different roles. Given $n>n_r,$ let us apply the induction
hypothesis to the partition $I_1,\ldots, I_{r-1}, \tilde{I}_r$ where the last set is now
$\tilde{I}_r = I_r\cup  \{n_{r}+1,\ldots, n\} = \{n_{r-1}+1,\ldots, n\}$ and the functions 
\begin{equation}
f_{n_r+1}=\ldots=f_{n}=0.
\label{zeros}
\end{equation} 
Let $\tilde{Z}^r$ denote (\ref{Zp}) corresponding to the set $\tilde{I}_r,$ i.e.
\begin{equation}
\tilde{Z}^r = \tilde{Z}^r(\sigma^1,\ldots,\sigma^{n}) =
\frac{ \exp \sum_{l\in \tilde{I}_r} F_l^r(\sigma^l,\sigma^1,\ldots,\sigma^{n})}
{\la\exp  F^r(\sigma,\sigma^1,\ldots,\sigma^{n})\ra_{\mathunderscore}^{|\tilde{I}_r|}}.
\label{tildeZr}
\end{equation}
Because of (\ref{zeros}),
$$
F^r(\sigma,\sigma^1,\ldots,\sigma^{n}) 
= 
f_1(\sigma\cdot\sigma^1)+\ldots+f_{n_r}(\sigma\cdot\sigma^{n_r})
=
F^r(\sigma,\sigma^1,\ldots,\sigma^{n_r})
$$
as in (\ref{Fp1}) for $p=r$, and for $l>n_r,$
$$
F_l^r(\sigma^l,\sigma^1,\ldots,\sigma^{n}) = F^r(\sigma^l,\sigma^1,\ldots,\sigma^{n_r}). 
$$
Therefore, the choice (\ref{zeros}) allows us to rewrite (\ref{tildeZr}) as
\begin{equation}
\tilde{Z}^r = Z^r(\sigma^1,\ldots,\sigma^{n_r}) 
\prod_{l=n_r+1}^{n}
\frac{\exp F^r(\sigma^l,\sigma^1,\ldots,\sigma^{n_r})}
{ \la \exp  F^r(\sigma,\sigma^1,\ldots,\sigma^{n_r}) \ra_{\mathunderscore} }.
\label{tildeZr2}
\end{equation}
This means that for fixed $\sigma^1,\ldots,\sigma^{n_r}$ the coordinates  $\sigma^l$ for $l>n_r$
are integrated with respect to the product measure
\begin{equation}
dG'(\sigma) = 
\frac{\exp F^r(\sigma,\sigma^1,\ldots,\sigma^{n_r})}
{ \la \exp  F^r(\sigma,\sigma^1,\ldots,\sigma^{n_r}) \ra_{\mathunderscore} }
dG(\sigma).
\label{Gprime}
\end{equation} 
If we denote by $\la\cdot\ra'$ the average over the coordinates $(\sigma^l)_{l>n_r}$ 
with respect to $(G')^{\otimes \infty}$  conditionally on $\sigma^1,\ldots, \sigma^{n_r}$
and if we choose $\Phi$ to be the function of the overlaps only on the coordinates
$(\sigma^l)_{n_r<l\leq n}$, then the induction hypothesis and (\ref{tildeZr2}) give
\begin{equation}
\e\la\Phi\ra = \e\la Z^1\ldots Z^r \la \Phi \ra' \ra.
\label{mainind}
\end{equation} 
We can think of $Z^1\ldots Z^r$ as the change of density and treat the functional of $\Phi$ 
on the right hand side as the probability on the overlaps of the coordinates  $(\sigma^l)_{l>n_r}$,
which again satisfies the Ghirlanda-Guerra identities.  We can now apply Theorem \ref{Th1} to
(or repeat its proof for) this functional to obtain the following. Let $n_{r+1}>n_r$ and
 $I_{r+1} = \{n_r+1,\ldots, n_{r+1}\}$. For $l\in I_{r+1}$ consider functions  $f_l(x) = f_l(R^{I_{r+1}},x)$,
define
\begin{equation}
F(\sigma,\sigma^{n_r+1},\ldots,\sigma^{n_{r+1}}) = 
f_{n_r+1}(\sigma\cdot\sigma^{n_r+1})+\ldots+f_{n_{r+1}} (\sigma\cdot\sigma^{n_{r+1}})
\label{F1ind}
\end{equation}
and for $l\in I_{r+1}$ define
\begin{equation}
F_l(\sigma,\sigma^{n_r+1},\ldots,\sigma^{n_{r+1}}) 
= 
F(\sigma,\sigma^{n_r+1},\ldots,\sigma^{n_{r+1}})
 - f_l(\sigma\cdot\sigma^l)+ \int \!\! f_l(x) \,d\mu(x).
\label{F2ind}
\end{equation}
Then, Theorem \ref{Th1} (or its proof) applied to the right hand side of (\ref{mainind}) implies
\begin{equation}
\e\la  \Phi \ra
=
\e\Bigl\la
Z^1\ldots Z^r 
\Bigl\la
\frac{\Phi \exp \sum_{l\in I_{r+1}} F_l(\sigma^l,\sigma^{n_r+1},\ldots,\sigma^{n_{r+1}}) }
{(\la\exp F(\sigma,\sigma^{n_r+1},\ldots,\sigma^{n_{r+1}}) \ra_{\mathunderscore}^\prime)^{|I_{r+1}|}}
\Bigr\ra'
\Bigr\ra,
\label{mainr}
\end{equation}
where $\la\cdot\ra^\prime_{\mathunderscore}$ is the average in $\sigma$ only with respect to the measure $G'$.
It remains to rewrite (\ref{mainr}) recalling the definition of the measure $G'$ in (\ref{Gprime}).
Since
\begin{align}
\la\exp F(\sigma,\sigma^{n_r+1},\ldots,\sigma^{n_{r+1}}) \ra_{\mathunderscore}^\prime
&=
\frac{\la \exp (F(\sigma,\sigma^{n_r+1},\ldots,\sigma^{n_{r+1}}) 
+ F^r(\sigma,\sigma^1,\ldots,\sigma^{n_r})) \ra_{\mathunderscore} }
{\la \exp  F^r(\sigma,\sigma^1,\ldots,\sigma^{n_r}) \ra_{\mathunderscore} }
\nonumber
\\
&=
\frac{\la \exp F^{r+1}(\sigma,\sigma^1,\ldots,\sigma^{n_{r+1}}) \ra_{\mathunderscore} }
{\la \exp  F^r(\sigma,\sigma^1,\ldots,\sigma^{n_r}) \ra_{\mathunderscore} }
\end{align}
by (\ref{Fp1}) and (\ref{F1ind}), the average $\la\cdot\ra'$ inside (\ref{mainr}) over the coordinates
$(\sigma^l)_{l\in I_{r+1}}$ can be rewritten as
\begin{align*}
&
\Bigl\la
\frac{\Phi \exp \sum_{l\in I_{r+1}} (F_l(\sigma^l,\sigma^{n_r+1},\ldots,\sigma^{n_{r+1}}) 
+F^r(\sigma^l,\sigma^{1},\ldots,\sigma^{n_{r}}) 
)}
{\la\exp F^{r+1}(\sigma,\sigma^{1},\ldots,\sigma^{n_{r+1}}) \ra_{\mathunderscore}^{|I_{r+1}|}}
\Bigr\ra
\\
&=
\Bigl\la
\frac{\Phi \exp \sum_{l\in I_{r+1}} F_l^{r+1}(\sigma^l,\sigma^{1},\ldots,\sigma^{n_{r+1}}) }
{\la\exp F^{r+1}(\sigma,\sigma^{1},\ldots,\sigma^{n_{r+1}}) \ra_{\mathunderscore}^{|I_{r+1}|}}
\Bigr\ra
= \la Z^{r+1} \Phi\ra,
\end{align*}
where we combined (\ref{F2ind}) and (\ref{Fp2}) in the numerator and where the average is taken only
over the coordinates $(\sigma^l)_{l\in I_{r+1}}$. This completes the induction step and finishes the proof.
\qed 

\noindent
It is worth formulating the general principle expressed in equations (\ref{Gprime}) and (\ref{mainind}) as a separate result.
Let $(\Omega,\Pr)$ be the probability space on which the random measure $G$ is defined and let $H$ be our Hilbert space.
Let ${\Pr}^{(n_r)}$ be the measure on $\Omega\times H^{n_r}$ defined by the change of density $Z^1\ldots Z^r$, i.e. for any
measurable sets $A\subseteq \Omega$ and $A_1,\ldots, A_{n_r}\subseteq H,$ probability
${\Pr}^{(n_r)} (A\times A_1\times \ldots\times A_{n_r})$ is given by 
\begin{equation}
\int_A \int_{A_1\times \ldots\times A_{n_r}}(Z^1\ldots Z^r)(\sigma^1,\ldots,\sigma^{n_r}) \,dG(\sigma^1)\ldots dG(\sigma^{n_r}) 
\,d\! \Pr(\omega).
\end{equation}
We can think of the random measure $G'$ in (\ref{Gprime}) as defined on $(\Omega\times H^{n_r},{\Pr}^{(n_r)})$ since it depends
on $\omega$ and $\sigma^1,\ldots, \sigma ^{n_r}$ and (\ref{mainind}) expresses the following.
\begin{theorem}\label{ThZr}
Under (\ref{GG}), the random measure $G'$ on $(\Omega\times H^{n_r},{\Pr}^{(n_r)})$ defined in (\ref{Gprime}) has the same 
distribution as the measure $G$ on $(\Omega,\Pr)$ in the sense that i.i.d. samples from these measures have the same joint 
overlap distributions.
\end{theorem}

\noindent
Let us write down another  generalization of Theorem \ref{Th1} in a slightly different direction on which our applications will be based.
Consider a finite index set $\A.$ Given $n\geq 1$ and configurations $\sigma^1,\ldots,\sigma^n,$ 
let $(B_\alpha)_{\alpha\in\A}$ be a partition of the Hilbert space $H$ such that for each $\alpha\in\A$ the indicator 
$I(\sigma\in B_\alpha)$ is a measurable function of $R^n$ and $(\sigma\cdot\sigma^l)_{l\leq n}$ and let
\begin{equation}
W_\alpha=W_\alpha(\sigma^1,\ldots,\sigma^n)=G(B_\alpha).
\label{WA}
\end{equation}
Let us define a map $T$ by
\begin{equation}
W=(W_\alpha)_{\alpha\in\A}\to T(W) = 
\Bigl(\frac{\la I_{B_\alpha} \exp F(\sigma,\sigma^1,\ldots,\sigma^n )\ra_{\mathunderscore}}
{\la \exp F(\sigma,\sigma^1,\ldots,\sigma^n )\ra_{\mathunderscore}} \Bigr)_{\alpha\in\A}.
\label{TA}
\end{equation}
We have the following invariance result for the weights $(W_\alpha)$ of the random partition $(B_\alpha)$.
\begin{theorem}\label{Th2}
Under (\ref{GG}), for any bounded measurable function $\varphi:\Reals^{n(n-1)/2}\times\Reals^{|\A|}\to \Reals$,
\begin{equation}
\e\la  \varphi(R^n, W)\ra
=
\e\Bigl\la
\frac{ \varphi(R^n,T(W)) \exp \sum_{l=1}^{n} F_l(\sigma^l,\sigma^1,\ldots,\sigma^n)}
{\la\exp F(\sigma,\sigma^1,\ldots,\sigma^n)\ra_{\mathunderscore}^n}
\Bigr\ra.
\label{nA}
\end{equation}
\end{theorem}\noindent
\textbf{Proof.} 
For each $\alpha\in\A$ let us take integer $n_\alpha\geq 0$ and let $m=n+\sum_{\alpha\in\A} n_\alpha.$ 
Let $(S_\alpha)_{\alpha\in\A}$ be any partition of $\{n+1,\ldots,m\}$ such that $|S_\alpha |=n_\alpha.$
Consider a continuous function $\Phi:\Reals^{n(n-1)/2}\to\Reals$ and let
$\Phi' = \Phi(R^n) \prod_{\alpha\in\A}\varphi_\alpha$ where 
$$
\varphi_\alpha = I (\sigma^l \in B_\alpha, \forall l\in S_\alpha )
$$
and let $f_l$ for $l\leq n$ be as in (\ref{F1}) and $f_{n+1}=\ldots=f_m=0.$ Let us now apply Theorem \ref{Th1} 
with these choices of functions $\Phi'$ and $f_l$ (and $n=m$). First of all, integrating out the coordinates 
$(\sigma^l)_{l>n}$, the left hand side of (\ref{main}) can be written  as
\begin{equation}
\e\la \Phi' \ra 
=
\e \bigl\la \Phi(R^n) \prod_{\alpha\in\A} \varphi_\alpha \bigr\ra
=
\e \bigl\la\Phi(R^n) \prod_{\alpha\in\A} W_\alpha^{n_\alpha}(\sigma^1,\ldots,\sigma^n) \bigr\ra
\label{lhscor}
\end{equation}
where $W_\alpha$'s were defined in (\ref{WA}). Let us now compute the right hand side of (\ref{main}). 
Since $f_{n+1}=\ldots=f_m=0,$ the denominator will be
$\bigl\la\exp F(\sigma,\sigma^1,\ldots,\sigma^n)\bigr\ra_{\mathunderscore}^m$  and
\begin{equation}
\sum_{l=1}^{m} F_l(\sigma^l,\sigma^1,\ldots,\sigma^m) 
= 
\sum_{l=1}^{n} F_l(\sigma^l,\sigma^1,\ldots,\sigma^n)
+
\sum_{l=n+1}^{m} F(\sigma^l,\sigma^1,\ldots,\sigma^n).  
\label{numcor}
\end{equation}
Since the denominator does not depend on  $(\sigma^{l})_{l>n}$,
integrating in the coordinate $\sigma^l$ for $l\in S_\alpha$ will produce a factor 
$$
\la I_{B_\alpha} \exp F(\sigma,\sigma^1,\ldots,\sigma^n )\ra_{\mathunderscore}.
$$
For each $\alpha\in \A$ we have $|S_\alpha| = n_\alpha$ of such coordinates and, therefore, the right hand side of (\ref{main}) 
is equal to
\begin{equation}
\e\Bigl\la
\frac{\Phi(R^n) \exp \sum_{l=1}^{n} F_l(\sigma^l,\sigma^1,\ldots,\sigma^n)}
{\la\exp F(\sigma,\sigma^1,\ldots,\sigma^n)\ra_{\mathunderscore}^n}
\prod_{\alpha\in\A} \Bigl(
\frac{\la I_{B_\alpha} \exp F(\sigma,\sigma^1,\ldots,\sigma^n )\ra_{\mathunderscore}}
{\la \exp F(\sigma,\sigma^1,\ldots,\sigma^n )\ra_{\mathunderscore}}
\Bigr)^{n_\alpha}
\Bigr\ra.
\label{rhscor}
\end{equation}
Comparing with (\ref{lhscor}), recalling (\ref{TA}) and approximating a continuous function $\phi$ on $[0,1]^{|\A|}$ by polynomials 
we get (\ref{nA}) first for products $\Phi(R^n) \phi(W)$, then for continuous functions $\varphi(R^n,W)$ and then for arbitrary
bounded measurable functions.
\qed

\section{Applications.}\label{Sec3}

Let us begin with the following special case of Theorem \ref{Th2}.
Let us consider some  sets $B_l$ for $l\leq n$. Let $\A$ be the power set of $\{1,\ldots,n\}$
and for each $\alpha\subseteq \{1,\ldots,n\}$ we define a set
\begin{equation}
B_\alpha(\sigma^1,\ldots,\sigma^n)= \{\sigma : \sigma\cdot\sigma^l \not\in B_l \Leftrightarrow l\in \alpha\}
\label{AI}
\end{equation}
which depends only on the overlaps. Given $t=(t_1,\ldots,t_n)\in\Reals^n$, let us now make the choice of functions
$f_l(R^n,x) = t_l I(x\not\in B_l).$ Since in this case, using notation $t_\alpha = \sum_{l\in \alpha} t_l$,
\begin{equation}
F(\sigma,\sigma^1,\ldots,\sigma^n)
=\sum_{l\leq n} t_l I(\sigma\cdot \sigma^l \not\in B_l) 
= \sum_{\alpha\in\A} t_\alpha I(\sigma \in B_\alpha)
\label{Fspec}
\end{equation}
we get
\begin{equation}
\bigl\la I_{B_\alpha} \exp F(\sigma,\sigma^1,\ldots,\sigma^n)\bigr\ra_{\mathunderscore}
=
W_\alpha e^{t_\alpha}
\end{equation}
and
\begin{equation} 
\bigl\la\exp F(\sigma,\sigma^1,\ldots,\sigma^n)\bigr\ra_{\mathunderscore}
=
\sum_{\alpha\in\A} W_\alpha e^{t_\alpha}
\end{equation}
If we denote $\Delta_t = \sum_{\alpha\in\A} W_\alpha e^{t_\alpha}$ the map $T$ in (\ref{TA}) becomes
\begin{equation}
W=(W_\alpha)_{\alpha\subseteq \{1,\ldots, n\}}\to T_t(W) = 
\Bigl(\frac{W_\alpha e^{t_\alpha}}{\Delta_t} \Bigr)_{\alpha\subseteq \{1,\ldots, n\}}.
\label{TI}
\end{equation}
Given a measurable function $\phi:\Reals^{2^n}\to \Reals$ and an arbitrary subset
\begin{equation}
B\subseteq \prod_{l< l'}^n B_l\cap B_{l'}\subseteq \Reals^{n(n-1)/2},
\label{setA}
\end{equation}
take the function $\varphi$ in (\ref{nA}) to be $\varphi = I(R^n\in B) \phi(W)$.
Since for $R^n\in B,$ the overlap $R_{l,l'}\in B_l\cap B_{l'}$ and therefore $f_l(R^n,R_{l,l'}) = 0$, we have
\begin{equation}
\sum_{l=1}^{n} F_l(\sigma^l,\sigma^1,\ldots,\sigma^n)
=
 \sum_{l\leq n} t_l \int \! I(x\not\in B_l)\,d\mu
=
 \sum_{l\leq n} t_l \,\mu((B_l)^c)
=:
\gamma_t.
\label{gammaI}
\end{equation}
With these choices of parameters, Theorem \ref{Th2} becomes:
\begin{theorem}\label{Th2a}
Under (\ref{GG}), for any bounded measurable function $\phi:\Reals^{2^n}\to \Reals$,
\begin{equation}
\e \la I(R^n\in B) \phi(W)\ra
=
\e\Bigl\la \frac{I(R^n\in B) \phi(T_t(W))  e^{\gamma_t}}{\Delta_t^n} \Bigr\ra.
\label{nI}
\end{equation}
\end{theorem}\noindent
To show how this implies Theorem \ref{Prop1}, let us consider the following special case with $n=2$.
Consider $q_2\leq q_1 \leq q\leq 1$ and let $B_l = [q_l,1]$ for $l\leq 2$ and $B=[q,1]$. 
The partition $(B_\alpha)$ will now consist of four sets
\begin{eqnarray}
&&
A_1(\sigma^1,\sigma^2) = \{\sigma : \sigma\cdot \sigma^1 \geq q_1, \sigma\cdot \sigma^2 < q_2\},
\nonumber
\\
&&
A_2(\sigma^1,\sigma^2) = \{\sigma : \sigma\cdot \sigma^1 < q_1, \sigma\cdot \sigma^2 \geq q_2\},
\nonumber
\\
&&
A_3(\sigma^1,\sigma^2) = \{\sigma : \sigma\cdot \sigma^1 \geq q_1, \sigma\cdot \sigma^2 \geq q_2\},
\nonumber
\\
&&
A_4(\sigma^1,\sigma^2) = \{\sigma : \sigma\cdot \sigma^1 < q_1, \sigma\cdot \sigma^2 < q_2\}.
\label{setsA}
\end{eqnarray}
Let $W_j = W_j(\sigma^1,\sigma^2)=G(A_j(\sigma^1,\sigma^2))$. With these notations,
\begin{equation}
\Delta_t = W_1 e^{t_2} + W_2 e^{t_1}+W_3+W_4 e^{t_1+t_2},
\label{Delta}
\end{equation}
the map
\begin{equation}
T_t(W) = \Bigl(\frac{W_1 e^{t_2}}{\Delta_t}, \frac{W_2 e^{t_1}}{\Delta_t}, 
\frac{W_3}{\Delta_t}, \frac{W_4 e^{t_1+t_2}}{\Delta_t}\Bigr)
\label{T}
\end{equation}
and $\gamma_t$ in (\ref{gammaI}) is
\begin{equation}
\gamma_t = t_1\p(R_{1,2}< q_1)+t_2\p(R_{1,2}< q_2).
\label{gamma}
\end{equation}
Theorem \ref{Th2a} gives that for any bounded measurable function $\phi:\Reals^4\to \Reals$,
\begin{equation}
\e \la I(R_{1,2}\geq q) \phi(W)\ra
=
\e\Bigl\la \frac{I(R_{1,2}\geq q) \phi(T_t(W))  e^{\gamma_t}}{\Delta_t^n} \Bigr\ra.
\label{n2}
\end{equation}
This readily implies Theorem \ref{Prop1}.

\smallskip
\noindent
\textbf{Proof of Theorem \ref{Prop1}.}
Let $s>0.$ Let us use (\ref{n2}) with $q_1=q_2=q,$ the choice of function
$$
\phi(W)=\frac{1}{(W_1 e^{-s} + W_2 e^{-s} +W_3 + W_4)^2}
$$
and the choice of $t_1=s, t_2=-s$. Then (\ref{n2}) becomes
$$
\e \Bigl\la \frac{I(R_{1,2}\geq q)}{(W_1 e^{-s} + W_2 e^{-s} +W_3 + W_4)^2} \Bigr\ra
=
\e \Bigl\la \frac{I(R_{1,2}\geq q)}{(W_1 e^{-2s} + W_2 +W_3 + W_4)^2} \Bigr\ra.
$$
Letting $s\to\infty$, by monotone convergence theorem we get
$$
\e \Bigl\la \frac{I(R_{1,2}\geq q)}{(W_3 + W_4)^2} \Bigr\ra
=
\e \Bigl\la \frac{I(R_{1,2}\geq q)}{(W_2 +W_3 + W_4)^2} \Bigr\ra.
$$
and condition (\ref{integ}) means that this quantity is finite. In that case, almost surely over
the choice of random measure $G$ and the choice of $(\sigma^1,\sigma^2)$ from $G^{\otimes 2}$,
if $\sigma^1\cdot\sigma^2 \geq q$ then $W_2=G(A_2(\sigma^1,\sigma^2)) = 0$ and by symmetry 
$W_1=0,$ which is another way to express  (\ref{ultra}).
\qed

\smallskip
\noindent
\textbf{Proof of Theorem \ref{Prop2}.}
In the setting of Theorem \ref{Th2}, let us take a partition consisting of three sets
$$
B_1 = \bigl\{\sigma : \bar{F}(\sigma,\sigma^1,\ldots,\sigma^n) > 0 \bigr\},\,
B_2 = \bigl\{\sigma : \bar{F}(\sigma,\sigma^1,\ldots,\sigma^n) < 0 \bigr\}
$$
and $B_3 = (B_1\cup B_2)^c$. Let $W_l=G(B_l)$ for $l \leq 3$. Let us apply (\ref{nA}) to the
function $\varphi(R^n, W) = I(W_1=0)$ and let us replace functions $f_l$ by $sf_l$ for $s>0.$ 
Since $(T(W))_1 = 0$ if and only if $W_1=0,$ we have $\varphi(R^n, T(W)) = I(W_1=0)$ 
and (\ref{nA}) becomes 
\begin{equation}
\e\la I(W_1=0)\ra
=
\e\Bigl\la
\frac{I(W_1=0) }
{\la\exp s\bar{F}(\sigma,\sigma^1,\ldots,\sigma^n)\ra_{\mathunderscore}^n}
\Bigr\ra.
\label{w1zero2}
\end{equation}
Whenever $W_1 = G(B_1) = 0,$ we have
$$
\la\exp sF(\sigma,\sigma^1,\ldots,\sigma^n)\ra_{\mathunderscore}
=
\la I_{B_2} \exp s\bar{F}(\sigma,\sigma^1,\ldots,\sigma^n)\ra_{\mathunderscore}+W_3
$$
and since $\bar{F}<0$ on $B_2,$ letting $s\to\infty$ implies by monotone convergence theorem
\begin{equation}
\e\la I(W_1=0)\ra
=
\e\Bigl\la
\frac{I(W_1=0)}{W_3^n}
\Bigr\ra.
\label{w1zero3}
\end{equation}
This shows that whenever $W_1=0$ we must have $W_3=1$ and,  therefore, $W_2=0.$ Similarly, one can show 
that whenever $W_2=0$ we must have $W_1=0$ which means that either both $W_1$ and $W_2$ are positive
or equal to zero. Of course, this implies that the measures $W_1+W_3$ 
of the set $\bar{F}\geq 0$ and $W_2+W_3$ of the set $\bar{F}\leq 0$ are always positive. 
\qed

\noindent
Let us give a partial generalization of Theorem \ref{Prop2} in the setting of Theorem \ref{Th1prime}.
Recall the notations of Theorem \ref{Th1prime} and define
\begin{equation}
\tilde{F}((\rho^l), (\sigma^l)) =
\sum_{p\leq r} \sum_{l\in I_p} F^p(\rho^l,\sigma^1,\ldots,\sigma^{n_p})
-
\sum_{p\leq r} \sum_{l\in I_p} F_l^p(\sigma^l,\sigma^1,\ldots,\sigma^{n_p}).
\end{equation}
The following holds.
\begin{theorem}\label{Prop2gen}
Under (\ref{GG}), with probability one over the choice of measure $G$, for  $G^{\otimes n_r}$-almost all  
$(\sigma^1,\ldots,\sigma^{n_r})$ the inequality $\tilde{F}\geq 0$ holds on a set of $(\rho^1,\ldots,\rho^{n_r})$ 
of positive measure $G^{\otimes n_r}$.
\end{theorem}\noindent
\textbf{Proof.}
Let us recall (\ref{Zprod}) and apply Theorem \ref{Th1prime} to the function $\Phi = 1$
and with functions $f_l$ replaced by $sf_l,$
\begin{equation}
1= \e\Bigl\la \frac{1} {\la \exp s \tilde{F}((\rho^l), (\sigma^l))  \ra_{\mathunderscore}}
\Bigr\ra. 
\end{equation}
If we consider the set
$$
B=\bigl\{(\sigma^1,\ldots, \sigma^{n_r}) : G\bigl(\{ (\rho^1,\ldots, \rho^{n_r}) : \tilde{F}((\rho^l), (\sigma^l)) \geq 0\}\bigr) = 0\bigr\}
$$
then
$$
1\geq \e\Bigl\la
\frac{I_B} {\la \exp s \tilde{F}((\rho^l), (\sigma^l)) \ra_{\mathunderscore}}
\Bigr\ra
$$
and, since on $B$ the limit $\lim_{s\to\infty}\la \exp s \tilde{F} \ra_{\mathunderscore} = 0$,
we must have $\e\la I_B\ra = 0$. This means that for $G^{\otimes n_r}$-almost all  $(\sigma^1,\ldots,\sigma^{n_r})$ 
the condition $\tilde{F} \geq 0$ is satisfied on the set of $(\rho^1,\ldots,\rho^{n_r})$ of positive measure $G^{\otimes n_r}$. 
\qed

\medskip
\noindent \textbf{Exchangeability and independence.} We will now prove Theorem \ref{ThEI}.
First, recall several properties of measures $G$ satisfying the Ghirlanda-Guerra identities (\ref{GG})
that will be used in the proof of Theorem \ref{ThEI}. We already mentioned in the introduction that if $q^*$ 
is the largest point in the support of the distribution of $\sigma^1\cdot\sigma^2$ under $\e G^{\otimes 2}$ 
then (\ref{GG}) implies that $G(\|\sigma\|^2 = q^*)=1$ with probability one which in the case of discrete 
measure $G=\sum_{l\geq 1} v_l \delta_{\xi_l}$ means that all $\|\xi_l\|^2 = q^*$ and, in particular, 
$R_{1,2} = \sigma^1\cdot\sigma^2 = q^*$ if and only if $\sigma^1 = \sigma^2.$ Also, by a well known
result of Talagrand (Section 1.2 in \cite{SG} or  Proposition  15.2.4 in \cite{SG2}) the weights $(v_l)$
 must have the Poisson-Dirichlet distribution $PD(\zeta)$ with $\zeta \in (0,1)$.
We recall that, given $\zeta\in (0,1),$ if $(u_{l})_{l \geq 1}$ is the decreasing enumeration of a Poisson 
point process on  $(0,\infty)$ with intensity measure $x^{-1- \zeta}dx$ on $(0,\infty)$ and $v_l = u_l/\sum_j u_j$ 
then the distribution of the sequence $(v_l)$ is called the Poisson-Dirichlet  distribution $PD(\zeta)$ (see e.g. \cite{PY}). 
We will be using the fact that 
\begin{equation}
\e\la I(\sigma^1\cdot\sigma^2 =1)\ra = \e \sum_{l\geq 1} v_l^2 = 1-\zeta
\label{zeta}
\end{equation}
(see e.g. \cite{Ruelle} or Section 13.1 in \cite{SG2}). 
When we sample $n$ replicas $\sigma^1,\ldots, \sigma^n$ from $G$, some of them could be equal so we can divide all indices
$\{1,\ldots, n\} = C_1 \cup \ldots \cup C_k$ into $k$ groups $C_1,\ldots, C_k$ such that $\sigma^l = \sigma^{l'}$ if and only if
$l$ and $l'$ belong to the same element of the partition. Let us for a moment fix one such partition and call it $C.$ We will use
the same notation to define the event
\begin{equation}
C = \bigl\{\forall l\not = l' \leq n,\,  R_{l,l'} = q^* \Longleftrightarrow \exists j\leq k \mbox{ such that } l,l'\in C_j \bigr \}.
\label{partitionC}
\end{equation}
Let $J=(j_1,\ldots, j_k)$ be a vector such that $j_l$ is the smallest index in $C_l$ and let us define 
\begin{equation}
R = (\sigma^{j_l} \cdot \sigma^{j_{l'}})_{l,l' \leq k} \,\,\mbox{ and }\,\, W = (G(\sigma^{j_1}),\ldots, G(\sigma^{j_k})).
\label{RW}
\end{equation}
Let us define the conditional distribution of $R$ and $W$ on $C$ by
\begin{equation}
\p_C(R\in A, W\in B) = \frac{\e\la I(R\in A, W\in B) I_C\ra}{\e\la I_C\ra }.
\label{PC}
\end{equation}
We will prove the following result from which Theorem \ref{ThEI} will easily follow.
\begin{theorem}\label{ThPC}
We have, 
\begin{equation}
\p_C(R\in A, W\in B) = \p_C(R\in A)\, \p_C(W\in B).
\label{indep}
\end{equation}
\end{theorem}\noindent
The main idea of the proof is contained in the computation in Lemma \ref{ThBeps} below which is
based on the invariance principle of Theorem \ref{Th2}.  Let 
\begin{equation}
{\cal W}_k = \{(w_1,\ldots, w_k) : \sum_{l\leq k} w_l < 1, w_1,\ldots, w_k>0\}.
\label{Wk}
\end{equation}
Given a vector $a=(a_1,\ldots, a_k)\in\Reals^k$ let us define $T_a: {\cal W}_k  \to {\cal W}_k$ by 
\begin{equation}
T_a(w) = \Bigl( \frac{w_1 e^{a_1}}{\Delta_a(w)},\ldots, \frac{w_k e^{a_k}}{\Delta_a(w)}\Bigr)
\label{Taw}
\end{equation}
where 
\begin{equation}
\Delta_a(w) = \sum_{l\leq k} w_l e^{a_l} + 1-\sum_{l\leq k} w_l. 
\end{equation}
One can easily check that for $a,b\in \Reals^k$ we have $T_a\circ T_b = T_{a+b}$ and, therefore, $T_a^{-1} = T_{-a}$.
Also, it is easy to check that 
\begin{equation}
\Delta_a(T_{-a}(w)) = \Delta_{-a}(w)^{-1}.
\label{DeltaT}
\end{equation}
Let us denote by $B_\eps(w)$ an open ball of radius $\eps$ centered at $w.$ Then the following holds.
\begin{lemma}\label{ThBeps}
For any $a=(a_1,\ldots, a_k)\in\Reals^k$ and $w\in {\cal W}_k$,
\begin{equation}
\lim_{\eps\to 0^+}  
\frac{\p_C(R\in A, W\in B_\eps(w))}{ \p_C(W\in B_\eps(w))}
=
\lim_{\eps\to 0^+}
 \frac{\p_C(R\in A, T_a(W) \in B_\eps(w))}{ \p_C(T_a(W) \in B_\eps(w))}
\label{Beps}
\end{equation}
whenever either of the limits exists.
\end{lemma}\noindent
\textbf{Proof.} For simplicity of notation, let us assume that $J=(j_1,\ldots, j_k) = \{1,\ldots, k\}$.
In (\ref{WA}) - (\ref{nA}) let us take $\A = \{1,\ldots, k+1\},$ $B_l = \{\sigma^{l}\}$ for $l\leq k$
and $B_{k+1} = \{\sigma^{1},\ldots,\sigma^{k}\}^c$, $f_{l}(x) = a_l I(x=q^*)$ for $l\leq k$ and
$f_l = 0$ for $l>k$. Then our notation $W_l = G(B_l) = G(\sigma^l)$ in (\ref{RW}) agrees with (\ref{WA})
for $l\leq k$ and we will forget about $W_{k+1}=G(B_{k+1})$ and only look at functions of 
$W=(W_1,\ldots, W_k).$ Let us take $\varphi(R^n, W)$ in (\ref{nA}) to be
$$
\varphi(R^n, W) = I(R\in A, W\in B) I(R^n\in C)
$$ 
where as in (\ref{RW}), $R = (\sigma^{l} \cdot \sigma^{{l'}})_{l,l' \leq k}$ is a $k\times k$ block in $R^n$.  
It is easy to check that with the choices we made, on the event $C$ the terms that appear on the right
hand side of (\ref{nA}) will become (recall (\ref{zeta}))
$$
\sum_{l=1}^{n} F_l(\sigma^l,\sigma^1,\ldots,\sigma^n) = \e\la I(R_{1,2} = q^*) \ra \sum_{l\leq k}a_l 
= (1-\zeta)  \sum_{l\leq k}a_l, 
$$
$$
\la \exp F(\sigma,\sigma^1,\ldots,\sigma^n) \ra_{\mathunderscore}
= \sum_{l\leq k} W_l e^{a_l} + 1-\sum_{l\leq k} W_l = \Delta_a(W)
$$
and the first $k$ coordinates of the map $T(W)$ in (\ref{TA}) are given by $T_a(W).$ Therefore, (\ref{nA}) implies
\begin{equation}
\e\la I(R\in A, W\in B) I_C\ra = \e \la I(R\in A, T_a(W)\in B) I_C Z_a(W) \ra
\end{equation}
where 
$$
Z_a(W) = \Delta_a(W)^{-n} \exp (1-\zeta)  \sum_{l\leq k}a_l
$$ 
and using this for $B= B_\eps(w)$,
\begin{equation}
\frac{\p_C(R\in A, W\in B_\eps(w))}{ \p_C(W\in B_\eps(w))}
=
\frac{\e\la I(R\in A, T_a(W) \in B_\eps(w)) I_C Z_a(W)\ra}{ \e\la I(T_a(W) \in B_\eps(w))I_C Z_a(W)\ra},
\label{BepsZ}
\end{equation}
assuming that the numerator is not zero. By (\ref{DeltaT}) and the fact that  $T_a^{-1} = T_{-a}$ we get that if 
$T_a(W) \in B_\eps(w)$ then $Z_a(W)$ takes values in the set  
$$
\Bigl\{\Delta_{-a}^n(w') \exp (1-\zeta)  \sum_{l\leq k}a_l : w'\in B_{\eps}(w)\Bigr\}.
$$
Therefore, as $\eps \to 0^+,$ $Z_a(W)$ converges uniformly over such $W$ to 
$\Delta_{-a}^n(w) \exp (1-\zeta)  \sum_{l\leq k}a_l$, a constant, which will cancel out on the right hand side
of (\ref{BepsZ}) and yield (\ref{Beps}).
\qed

\noindent \textbf{Proof of Theorem \ref{ThPC}.}
Since the weights $(v_l)$ have Poisson-Dirichlet distribution $PD(\zeta)$, using well known 
representations of the Poisson point process with intensity measure $x^{-1-\zeta} dx$ and the corresponding
representation of $(v_l)$ (see, e.g. Proposition 8 in \cite{PY}), one can easily check that the distribution of 
any finite subset of $k$ weights is absolutely continuous with respect to the Lebesgue measure on $\Reals^k$.
This implies that the distribution of $W$ in (\ref{RW})  under $\p_C$ is also absolutely continuous with respect 
to the Lebesgue measure on $\Reals^k$ since, on the event $C,$ $\sigma^{j_l}$ are all different for $l\leq k$.
Let $p(w)$ be the Lebesgue density of this distribution and let $p_A(w)$ be the conditional expectation 
of $I(R\in A)$ given $W$ under $\p_C.$ Then, for any measurable set $B$ on $\Reals^k$,
\begin{equation}
\p_C(R\in A, W\in B) = \int_B p_A(w) p(w) dw,\,\,
\p_C(W\in B) = \int_B p(w) dw.
\label{pAp}
\end{equation}
To prove (\ref{indep}), it is enough to show that $p_A(w)$ is a constant a.e. on the set $\{w: p(w)>0\}$.
By the Lebesgue differentiation theorem, for almost every $w'\in\Reals^k$ one has (Corollary 1.6 in \cite{Stein})
\begin{equation}
\lim_{\eps\to0^+} \frac{1}{|B_\eps(w')|}\int_{B_\eps(w')} |p_A(w)p(w) - p_A(w')p(w')|dw = 0
\label{Leb1}
\end{equation}
and
\begin{equation}
\lim_{\eps\to0^+} \frac{1}{|B_\eps(w')|}\int_{B_\eps(w')} |p(w) - p(w')|dw = 0.
\label{Leb2}
\end{equation}
If $p_A(w)$ is not a constant a.e. on $\{p(w)>0\}$ then we can find two points $w',w''$ for which
both (\ref{Leb1}) and (\ref{Leb2}) hold and such that $p(w'), p(w'')>0$ and  $p_A(w')\not = p_A(w'').$
We can also assume that $w',w''\in {\cal W}_k$ in (\ref{Wk}) since $\p_C(W\not\in {\cal W}_k) = 0.$ First of all, equations
(\ref{pAp}) - (\ref{Leb2}) imply that the left hand side of (\ref{Beps})
\begin{equation}
\lim_{\eps\to 0^+}  
\frac{\p_C(R\in A, W\in B_\eps(w'))}{ \p_C(W\in B_\eps(w'))}
= p_A(w').
\label{BepsLHS}
\end{equation}
It is easy to check that if we take
$$
a_l = \log \frac{w_l'}{w_l''} - \log \frac{1-w_1'-\ldots - w_k'}{1- w_1''-\ldots - w_k''}
$$
for $l\leq k$ then $T_a(w'') = w'$ for $T_a$ defined in (\ref{Taw}). Equations (\ref{Beps}) and (\ref{BepsLHS}) imply that
\begin{equation}
\lim_{\eps\to 0^+}   \frac{\p_C(R\in A, W\in T_{-a} B_\eps(w'))}{ \p_C(W\in T_{-a} B_\eps(w'))} = p_A(w').
\label{contr1}
\end{equation}
To finish the proof, we will follow the argument of Corollary 1.7 in \cite{Stein} and use the fact 
that the sets $T_{-a} (B_\eps(w'))$ are of bounded eccentricity.
Since all partial derivatives of $T_a$ are uniformly bounded in a small neighborhood of $w''$
and all partial derivatives of $T_a^{-1}=T_{-a}$ are uniformly bounded in a small neighborhood of $w'$, there
exist constants $c, C>0$ such that $B_{c \eps}(w'') \subseteq T_{-a} (B_\eps(w')) \subseteq B_{C \eps}(w'')$ 
for small $\eps>0$. Therefore,
\begin{align*}
\frac{1}{| T_{-a} (B_\eps(w')) |}\int_{ T_{-a} (B_\eps(w'))} |p(w) - p(w'')|dw
& \leq
\frac{1}{ |B_{c\eps}(w'')|}\int_{B_{C\eps}(w'')} |p(w) - p(w'')|dw 
\\
& = \frac{C^k}{c^k |B_{C\eps}(w'')|}\int_{B_{C\eps}(w'')} |p(w) - p(w'')|dw
\end{align*}
and using that (\ref{Leb2}) holds for $w''$ implies
$$
\lim_{\eps\to0^+} \frac{1}{| T_{-a} (B_\eps(w')) |}\int_{ T_{-a} (B_\eps(w'))} |p(w) - p(w'')|dw = 0.
$$
Similarly, using (\ref{Leb1}) for $w''$ we get
$$
\lim_{\eps\to0^+} \frac{1}{| T_{-a} (B_\eps(w')) |}\int_{ T_{-a} (B_\eps(w'))} |p_A(w) p(w) - p_A(w'') p(w'')|dw = 0.
$$
These equations together with (\ref{pAp}) for $B = T_{-a} (B_\eps(w'))$ imply that
$$
\lim_{\eps\to 0^+}   \frac{\p_C(R\in A, W\in T_{-a} B_\eps(w'))}{ \p_C(W\in T_{-a} B_\eps(w'))} = p_A(w'')
$$
and, recalling (\ref{contr1}), we get $p_A(w') = p_A(w'')$ - a contradiction.
\qed

\noindent \textbf{Proof of Theorem \ref{ThEI}}.
Let us fix $m\geq 1.$ Let $\pi = (\pi_1,\ldots, \pi_m)$ denote the vector of indices corresponding to the $m$ largest
{\it different} weights among $G(\sigma^1),\ldots, G(\sigma^n).$ This assumes that the partition $C$ defined before
(\ref{partitionC}) has at least $m$ elements, so let us denote by $\C_m$ all such partitions. Recalling the vector
$J$ in (\ref{RW}), we can assume that if $\pi=I=(i_1,\ldots,i_m)$ then $i_1,\ldots, i_m \in \{j_1,\ldots, j_k\}$ and 
let us denote this fact by writing $I \subset J(C)$, where we made the dependence of $J=J(C)$ on $C$ explicit.
Let us now also make the dependence of $R=R^J$ and $W=W^J$ in (\ref{RW}) on $J$ explicit and let 
$R^I = (\sigma^{i_l}\cdot \sigma^{i_{l'}})_{l,l'\leq m}$ and $W^I = (G(\sigma^{i_l}))_{l\leq m}$. 
Given a subset $B \subseteq {\cal W}_m $, we can rewrite the event $\{\pi=I, W^I\in B\}$ in terms of $W^J$
by defining a set $B_I \subseteq {\cal W}_k$ such that
$$
\{W^J \in B_I\} = \{W^I \in B, W^I \mbox{ is the vector of largest $m$ weights in }  W^J\}.
$$
Then we can write
\begin{align}
\e\la I(R^\pi \in A, W^\pi\in B)\ra
& =
\sum_{C\in\C_m} \sum_{I\subset J(C)} \e\la I(R^I\in A, W^I \in B,\pi=I) I_C\ra
\nonumber
\\
& =
\sum_{C\in\C_m} \sum_{I\subset J(C)} \e\la I(R^I\in A, W^J \in B_I) I_C\ra
\nonumber
\\
\mbox{ (by Theorem \ref{ThPC}) }
& = 
\sum_{C\in\C_m} \sum_{I\subset J(C)} \p_C (R^I\in A) \,\e\la I(W^J \in B_I) I_C\ra.
\label{sums}
\end{align}
We will show in a second that $\p_C(R^I \in A)$ depends on $C$ and $I$ only through $m$ and if we denote 
$D_m = \{\sigma^1,\ldots, \sigma^m \mbox{ - all different}\}$ and $R^m = (\sigma^l \cdot \sigma^{l'})_{l,l'\leq m}$ then
\begin{equation}
\p_C(R^I \in A) = \p_{D_m}(R^m\in A) := \frac{\e \la I(R^m\in A) I_{D_m}\ra}{\e\la I_{D_m}\ra}.
\label{pDm}
\end{equation}
Then, using (\ref{pDm}) in (\ref{sums}) we get
$$
\e\la I(R^\pi \in A, W^\pi\in B)\ra = \p_{D_m}(R^m \in A) \,\e\la I(W^\pi\in B)\ra .
$$
When $n$ gets large, with high probability the sample $\sigma^1,\ldots,\sigma^n$ from $G$ will
contain points $\xi_1,\ldots, \xi_m$ corresponding to the largest weights $v_1,\ldots, v_m$ in $G$.
Therefore,
\begin{align}
\p((\xi_l \cdot \xi_{l'})_{l,l'\leq m}\in A, (v_l)_{l\leq m}\in B)
&=
\lim_{n\to\infty} \e\la I(R^\pi \in A, W^\pi\in B)\ra
\nonumber
\\
&=
\p_{D_m}(R^m \in A) \lim_{n\to\infty} \e\la I(W^\pi\in B)\ra 
\nonumber
\\
& = 
 \p_{D_m}(R^m \in A) \, \p( (v_l)_{l\leq m}\in B). 
\end{align}
This proves that $(\xi_l \cdot \xi_{l'})_{l,l'\leq m}$ is independent of  $(v_l)_{l\leq m}$ and its distribution
is invariant under permutations of coordinates since the distribution of $R^m$ under $\p_{D_m}$ is 
obviously invariant under permutations of coordinates. It remains to explain why (\ref{pDm}) holds.
Again, for simplicity of notation, suppose that $I=\{1,\ldots, m\}$ and $J(C)=\{1,\ldots, k\}$. Suppose first that
$k<n$ and assume, without loss of generality, that
$n\in C_1$ - the element of the partition such that $1\in C_1$. Let us denote by $C_l' = C_l \cap \{1,\ldots, n-1\}$
and define the event 
\begin{equation}
C' = \bigl\{\forall l\not = l' \leq n-1,\,  R_{l,l'} = q^* \Longleftrightarrow \exists j\leq k \mbox{ such that } l,l'\in C_j'  \bigr \}.
\label{partitionC2}
\end{equation}
Then, clearly, $I_C = I_{C'} I(R_{1,n} = q^*)$ and (\ref{GG}) implies that
$$
\e \la I(R^I\in A) I_{C}\ra = \e \la I(R^I\in A) I_{C'} I(R_{1,n} = q^*)\ra
=
\frac{|C_1| - 1 -\zeta}{n-1} \e \la I(R^I\in A) I_{C'} \ra.
$$
Similarly,
$$
\e \la I_{C}\ra =  \frac{|C_1| - 1 -\zeta}{n-1} \e \la  I_{C'} \ra
$$
and, therefore, $\p_C(R^I \in A) = \p_{C'}(R^I \in A). $ We can continue to remove coordinates
outside of $J(C)$ one by one until $k=n$. Once we are left with different configurations $\sigma^1,\ldots,\sigma^k$
we can remove in a similar fashion coordinates with indices $m+1,\ldots, k$ to finish the proof of (\ref{pDm}),
which competes the proof of Theorem \ref{ThEI}.
\qed

\smallskip
\noindent
\textbf{Identities for the Poisson-Dirichlet distribution}.
Finally, let us write down a couple of straightforward consequences of Theorems \ref{Th1} and \ref{Th1prime} for 
the Poisson-Dirichlet distribution $PD(\zeta)$.  
Due to a result of Talagrand that we mentioned above (see e.g. Section 1.2 in \cite{SG} or Theorem 15.2.1 in \cite{SG2}), 
the simplest measure for which the Ghirlanda-Guerra identities hold is the discrete measure concentrated 
on the orthonormal basis $(e_k)$ with weights $v_k=G(\{e_k\})$ from the Poisson-Dirichlet distribution 
$PD(\zeta)$ for $\zeta\in(0,1).$  First, let us write down what Theorem \ref{Th1} 
says for this measure. Since the overlap now takes only two values $0$ and $1,$ any function of $n$ configurations
depends only on their partition into equal configurations. Therefore, we only need to write down what happens 
for any such particular partition. Given $r\geq 1$, let $I_1,\ldots,I_r$ be a partition of $\{1,\ldots,n\}$ and let 
$n_p = |I_p|$ for $p\leq r$. Let $\Phi(R^n)$ be the indicator of the set 
\begin{equation}
\bigl\{(\sigma^1,\ldots,\sigma^n) : \sigma^j = \sigma^{j'} \Leftrightarrow j,j'\in I_p \mbox{ for some } p\leq r\bigr\}.
\label{event}
\end{equation}
Since the overlap takes only two values, the most general choice of functions $f_j$ that we can make here is $f_j(x) = t_j I(x=1)$ 
for some $t_j \in\Reals.$  Since there is one-to-one correspondence between configurations in the set (\ref{event}) and $r$ different
indices $l_1\not = \ldots\not = l_r \in \Natural$ such that $\sigma^j = e_{l_p}$ for $j\in I_p$, we can rewrite (\ref{main}) in terms of 
$(v_{l_1},\ldots,v_{l_r})$. Using (\ref{zeta}) and letting $s_p = \sum_{j\in I_p} t_j,$ one can easily check that (\ref{main}) can be written as
\begin{equation}
\e\sum_{l_1\not = \ldots\not = l_r} v_{l_1}^{n_1}\ldots v_{l_r}^{n_r}  
=
\e\sum_{l_1\not = \ldots\not = l_r} 
\frac{e^{\sum_{p\leq r} (n_p-\zeta ) s_p}}
{(\sum_{p\leq r} v_{l_p} e^{s_p} + 1-\sum_{p\leq r}v_{l_p})^n}\,\,
v_{l_1}^{n_1}\ldots v_{l_r}^{n_r} .
\label{PDex}
\end{equation}
For example, when $n=2, r=2$ and $I_1=\{1\}, I_2=\{2\}$, $t_1=-t_2=t,$ (\ref{PDex}) becomes
$$
\e \sum_{l\not =l'} v_l v_{l'} = \e \sum_{l\not = l'} \frac{v_l v_{l'}}{(v_l e^t + v_{l'}e^{-t} + 1-v_l -v_{l'})^2}.
$$
Notice that one can not take the formal limit $t\to\infty$ on the right hand side for lack of integrability. 
To give another example,  if in the notations of Theorem \ref{Th1prime} we take $r=2,$ $I_1 = \{1\}, I_2=\{2\}, f_1=f_2 = t I(x<1)$ 
and $\Phi=1$ then (\ref{maingen}) becomes
\begin{eqnarray*}
1
&=&
\e\sum_{l\not = l'} \frac{ v_l v_{l'}  e^{2\zeta t}}{(v_l +e^t(1-v_l)(v_l+v_{l'} + e^t(1-v_l-v_{l'})))}
\\
&&
+\,\,\e\sum_{l\geq 1}  \frac{v_l^2 e^{2\zeta t}}{(v_l +e^t(1-v_l))(v_l+e^{2t}(1-v_l))}.
\end{eqnarray*}
Similarly to (\ref{PDex}), one can also write down the general case of Theorem \ref{Th1prime} for the Poisson-Dirichlet distribution, 
but we will omit the details here.
\qed

\end{document}